\documentstyle[txmac,a4,
amssymb,%
twoside,%
nocaphead,%
epsf,%
mypic,times,mathptm]{article}

\advance\oddsidemargin by -1.8cm
\advance\evensidemargin by -1.8cm
\advance\textwidth by 3.6cm

\def\mynewtheo#1#2{%
\newtheorem{@#1}{#2}
\newenvironment{#1}{\begin{@#1}\rm}{\end{@#1}}}

\mynewtheo{lemma}{Lemma}
\mynewtheo{exer}{Exercise}
\mynewtheo{theo}{Theorem}
\mynewtheo{rem}{Remark}
\mynewtheo{problem}{Problem}
\mynewtheo{defi}{Definition}
\mynewtheo{conj}{Conjecture}
\mynewtheo{corr}{Corollary}
\mynewtheo{prop}{Proposition}
\mynewtheo{question}{Question}
\mynewtheo{exam}{Example}
\mynewtheo{warn}{Warning}
\mynewtheo{conv}{Convention}

\makeatletter

\newenvironment{theorem}{\begin{theo}}{\end{theo}}

\newenvironment{eqn}{\begin{equation}}{\end{equation}}

\begin{document}

\makeatletter

\def\bysame{\same[\kern2cm]\,}
\def\lra{\longrightarrow}
\def\Ra{\Rightarrow}
\def\La{\Leftarrow}
\def\Md{\max\deg}
\def\md{\min\deg}
\def\qed{\hfill\@mt{\Box}}
\def\@mt#1{\ifmmode#1\else$#1$\fi}
\def\qqed{\hfill\@mt{\Box\enspace\Box}}
\def\TM{$^\text{\raisebox{-0.2em}{${}^\text{TM}$}}$}
\def\restr#1{\raisebox{-1mm}{$\big|_{#1}$}}

\let\ap\alpha
\let\nb\nabla
\let\tl\tilde
\let\lm\lambda
\let\Lm\Lambda
\let\dl\delta
\let\Dl\Delta
\let\eps\varepsilon
\let\gm\gamma
\let\bt\beta
\let\th\theta
\def\db{{\dot\beta}}
\let\Sg\Sigma
\let\sg\sigma
\def\Ev{{\operator@font Ev}\,}
\def\ASym{{\operator@font ASym}\,}
\def\Sym{{\operator@font Sym}\,}
\def\Aut{{\operator@font Aut}\,}
\def\dig{{\operator@font diag}\,}
\def\tr{{\operator@font tr}\,}
\def\Ker{{\operator@font Ker}\,}
\def\rk{{\operator@font rk}\,}
\def\cE{{\cal E}}
\def\cV{{\cal V}}
\def\cX{{\cal X}}
\def\cL{{\cal L}}
\def\cB{{\cal B}}
\def\cF{{\cal F}}
\def\bZ{{\Bbb Z}}
\def\bC{{\Bbb C}}
\def\bR{{\Bbb R}}
\def\bQ{{\Bbb Q}}
\def\bN{{\Bbb N}}
\def\hH{{\hat H}}
\def\fg{{\frak g}}
\def\tE{{\tl E}}
\def\gl{{\mathfrak gl}}
\def\sl{{\mathfrak sl}}
\let\ds\displaystyle
\let\vd\perp
\let\es\enspace
\let\sg\sigma
\let\dt\det
\let\op\oplus
\let\sm\setminus
\let\ol\overline
\def\epsfs#1#2{{\catcode`\_=11\relax\ifautoepsf\unitxsize#1\relax\else
\epsfxsize#1\relax\fi\epsffile{#2.eps}}}
\def\epsfsv#1#2{{\vcbox{\epsfs{#1}{#2}}}}
\def\vcbox#1{\setbox\@tempboxa=\hbox{#1}\parbox{\wd\@tempboxa}{\box
  \@tempboxa}}
\def\p{\epsfsv{2cm}}
\def\br#1{\left\lfloor#1\right\rfloor}
\def\ang#1{\left\langle#1\right\rangle}

\def\@test#1#2#3#4{%
  \let\@tempa\go@
  \@tempdima#1\relax\@tempdimb#3\@tempdima\relax\@tempdima#4\unitxsize
    \relax
  \ifdim \@tempdimb>\z@\relax
    \ifdim \@tempdimb<#2%
      \def\@tempa{\@test{#1}{#2}}%
    \fi
  \fi
  \@tempa
}

\def\go@#1\@end{}
\newdimen\unitxsize
\newif\ifautoepsf\autoepsftrue

\unitxsize4cm\relax
\def\epsfsize#1#2{\epsfxsize\relax\ifautoepsf
  {\@test{#1}{#2}{0.1 }{4   }
		{0.2 }{3   }
		{0.3 }{2   }
		{0.4 }{1.7 }
		{0.5 }{1.5 }
		{0.6 }{1.4 }
		{0.7 }{1.3 }
		{0.8 }{1.2 }
		{0.9 }{1.1 }
		{1.1 }{1.  }
		{1.2 }{0.9 }
		{1.4 }{0.8 }
		{1.6 }{0.75}
		{2.  }{0.7 }
		{2.25}{0.6 }
		{3   }{0.55}
		{5   }{0.5 }
		{10  }{0.33}
		{-1  }{0.25}\@end
		\ea}\ea\epsfxsize\the\@tempdima\relax
		\fi
		}

\let\old@tl\~
\def\~{\raisebox{-0.8ex}{\tt\old@tl{}}}

\author{A. Stoimenow\footnotemark[1]\\[2mm]
\small Department of Mathematical Sciences, \\
\small KAIST, Daejeon, 307-701, Korea \\
\small e-mail: {\tt stoimeno@mathsci.kaist.ac.kr}\\
\small WWW: {\hbox{\web|http://mathsci.kaist.ac.kr/~stoimeno/|}}
\\[3mm]
}

{\def\thefootnote{\fnsymbol{footnote}}
\footnotetext[1]{Supported by BK21 Project.}
\def\thefootnote{}
\footnotetext[2]{
\ry{1.3em}\em{AMS subject classification:} 57M25 (primary), 22E10,
32M05, 20F36, 11E39 (secondary)\\
\rx{1.8em}%
\em{Keywords:} link, Lie group, maximal subgroup, braid, 
representation.\vspace{1mm}
}}

\title{\large\bf \uppercase{The density of Lawrence-Krammer and}
\\[3mm]
\uppercase{non-conjugate braid representations of links}
\\[2mm]
}

\date{\phantom{%
\large Current version: \curv\ \ \ First version:
\makedate{25}{7}{2008}}\vspace{-8mm}
}

\newcounter{type}%
\def\labeltype{(\roman{type})}%
\def\mylabel{\hbox to 7mm{\refstepcounter{type}\labeltype
\xdef\@currentlabel{\labeltype}\hfill}}%

\maketitle

\long\def\@makecaption#1#2{%
   \vskip 10pt
   {\let\label\@gobble
   \let\ignorespaces\@empty
   \xdef\@tempt{#2}%
   }%
   \ea\@ifempty\ea{\@tempt}{%
   \setbox\@tempboxa\hbox{%
      \fignr#1#2}%
      }{%
   \setbox\@tempboxa\hbox{%
      {\fignr#1:}\capt\ #2}%
      }%
   \ifdim \wd\@tempboxa >\captionwidth {%
      \rightskip=\@captionmargin\leftskip=\@captionmargin
      \unhbox\@tempboxa\par}%
   \else
      \hbox to\captionwidth{\hfil\box\@tempboxa\hfil}%
   \fi}%
\def\fignr{\small\sffamily\bfseries}%
\def\capt{\small\sffamily}%

\newdimen\@captionmargin\@captionmargin2cm\relax
\newdimen\captionwidth\captionwidth\hsize\relax

\let\reference\ref
\def\eqref#1{(\protect\ref{#1})}

\def\proof{\@ifnextchar[{\@proof}{\@proof[\unskip]}}
\def\@proof[#1]{\noindent{\bf Proof #1.}\enspace}

\def\myfrac#1#2{\raisebox{0.2em}{\small$#1$}\!/\!\raisebox{-0.2em}
  {\small$#2$}}
\def\abstractname{}

\def\inx#1#2{\setbox\@tempboxa\hbox{$#1$}\parbox[t]{\wd\@tempboxa}%
{$#1$\\\scriptsize\hbox to \wd\@tempboxa{\hss $#2$\hss}}}

\@addtoreset {footnote}{page}

\renewcommand{\section}{%
   \@startsection
         {section}{1}{\z@}{-1.5ex \@plus -1ex \@minus -.2ex}%
               {1ex \@plus.2ex}{\large\bf}%
}
\renewcommand{\@seccntformat}[1]{\csname the#1\endcsname .
\quad}

{\let\@noitemerr\relax
\vskip-2.7em\kern0pt\begin{abstract}
\noindent{\bf Abstract.}\enspace
We use some Lie group theory and Budney's unitarization of
the Lawrence-Krammer representation, to prove that for generic
parameters of definite form the image of the representation (also
on certain types of subgroups) is dense
in the unitary group. This implies that, except possibly for closures
of full-twist braids, all links have infinitely many conjugacy
classes of braid representations on any non-minimal number of (and
at least 4) strands. 
\\[2mm]
\end{abstract}
}

\section{Introduction}



The Lawrence-Krammer representation \cite{Lawrence,Krammer,Bigelow2}
$\rho_n$
of the braid group $B_n$ into $SL(p,\bZ[t^{\pm 1}q^{\pm 1}])$, with
$p=n(n-1)/2$, has become recently of interest as the first faithful
representation of braid groups. In this paper, we are concerned with
the identification of the image of $\rho_n$ on $B_n$ and certain
types of subgroups thereof. An important
property, unitarizability, is found by Budney \cite{Budney} (see
theorem \reference{TB} below). Our main result is the following.

\begin{theorem}\label{t1}
Assume $q,t$ with $|t|=|q|=1$ are chosen so that $t^aq^b=1$ for
$a,b\in\bZ$ implies that $a=b=0$, and the Budney form is definite
at $q,t$. Moreover, assume that $\rho_n$ is irreducible at $q,t$.
Then $\ol{\rho_n(B_n)}= U(p)$ (for $p=n(n-1)/2$).
\end{theorem}

This is analogous to a previous result for the \em{Burau
representation} $\psi_n$ (which we simply call `Burau' below)
in \cite{reiko}.

The irreducibility of $\rho_n$ will be treated extra with lemma
\ref{t1_}. It should be pointed out that it has been proved at
separate places.
There is written account by M.~Zinno \cite{Zinno}, though it was
observed also by others, incl. V.~Jones, R.~Budney and W.~T.~Song.
W.~T.~Song has proved the stronger statement that Budney's form is
the only unitarizing form. (He informed me of this result prior to
my proof of theorem \reference{t1}, although now the theorem
implies this uniqueness property, at least for definite form,
by remark \reference{z2}.) But
it appears all this material is (yet) unpublished. There is closely
related work of I.~Marin \cite{Marin}, which we discuss in \S\ref{S6}.
The proof of lemma \reference{t1_} is provided for completeness and
because of its simplicity compared to other methods.

Our main motivation was again the study of braid representations
of links. The problem to determine conjugacy classes of braid
representations of a given link goes back to the 60s. For some
early work see e.g. \cite{MT}. With the increasing attention given
to braids the problem was studied later e.g. in \cite{Fiedler,%
Morton,Fukunaga}. We apply theorem \reference{t1} to prove

\begin{theorem}\label{t2}
Assume $L$ is a link and $n>b(L)$. Then there exist infinitely many
conjugacy classes of $n$-braid representations of $L$, except if\\
{\def\labeltype{(\alph{type})}%
\mylabel $n\le 3$ or (possibly) \label{t2:a}\\ 
\mylabel $L$ is a $(n-1,k(n-1))$-torus link $(k\in\bZ)$. (This includes
the case $k=0$ of the $n-1$-component trivial link.) \label{t2:b}\\
\setcounter{type}{0}%
}
\end{theorem}

The number $b(L)$ is the minimal number of strands of a braid
representation of $L$, and is called \em{braid index} (see
e.g. \cite{Murasugi}).
The case \ref{t2:a} is very well-known from \cite{BM2} to need to be
excluded, but we do not know anything about whether any link of case
\ref{t2:b} is indeed exceptional. Still the theorem almost completely
settles the (in)finiteness for $n>b(L)$. For $n=b(L)$ the situation
is far more complicated; there are certain links also for $n\ge 4$
with a single conjugacy class, e.g. unlinks \cite{BM}, and a further
example due to Ko and Lee. Contrarily, Shinjo has announced to me
(in preparation) that she has, as extension of her work \cite{Shinjo},
found very general families of knots with infinitely many minimal
braid conjugacy classes. It is possible that the decision problem
when finitely many and when infinitely many classes occur for $n=b(L)$
is too complex to have a meaningful answer.

%

Most of the rest of the paper, until the end of \S\reference{S5},
will be devoted to the proof of theorem \reference{t1}. The proof
is rather Lie-group theoretic, and we will need to bring up some
related material along the way. In \S\reference{S6} we extend
theorem \reference{t1} to denseness of the image of subgroups of
$B_n$. In \S\reference{S7} we discuss theorem \reference{t2}.

\section{Lawrence-Krammer representation and its unitarization}

The $n$-strand \em{braid group} $B_n$ is considered generated by
the Artin \em{standard generators} $\sg_i$ for $i=1,\dots,n-1$
\cite{Artin,Garside}. These are subject to relations of the
type $[\sg_i,\sg_j]=1$ for $|i-j|>1$, which we call
\em{commutativity relations} (the bracket denotes the
commutator) and $\sg_{i+1}\sg_i\sg_{i+1}= \sg_i\sg_{i+1}\sg_i$,
which we call \em{Yang-Baxter} (or shortly YB) \rm{relations}. We
write $[\bt]$ for the \em{exponent sum} of $\bt$, and set $B_{k,l}
\subset B_n$ for $1\le k<l\le n$ to be the subgroup $<\,\sg_k,
\dots,\sg_{l-1}\,>$ (where angle brackets mean `generated by').

The representation $\rho_n$ of $B_n$ can be defined as operating
on a complex vector space $R=\bC^p$ with $p=n(n-1)/2$ with
basis $\{\,v_{i,j}\,:\,1\le i<j\le n\,\}$ by
\begin{eqn}\label{rf}
\rho_n(\sigma_i) v_{j,k} = \left\{
\begin{array}{lr}
v_{j,k} & i\notin \{j-1,j,k-1,k\}, \\
qv_{i,k} + (q^2-q)v_{i,j} + (1-q)v_{j,k} & i=j-1, \\
v_{j+1,k} & i=j\neq k-1, \\
qv_{j,i} + (1-q)v_{j,k} - (q^2-q)tv_{i,k} & i=k-1\neq j,\\
v_{j,k+1} & i=k,\\
-tq^2v_{j,k} & i=j=k-1.
\end{array}
\right.
\end{eqn}
Here $t$ and $q$ may \em{a priori} be arbitrary non-zero
complex numbers. However, we will choose them always so
that $|t|=|q|=1$. (We will sometimes write $q,t$ explicitly
as parameters of $\rho_n$, with the understanding that
a braid cannot be confused with a complex number.) The
reason is the following result, which is of main importance
below.

\begin{theorem}(Budney \cite{Budney})\label{TB}
The Lawrence-Krammer representation unitarizable if $|q|=|t|=1$. 
\end{theorem}

In other words, for
such $t$ and $q$, Budney \cite{Budney} defines a \em{unitarizing}
form $<\,.\,,\,.\,>$ of $\rho_n$ on $\bC^p$. This is a sesquilinear
pairing respected by the action of $\rho_n$, in the sense that for all
$\bt\in B_n$ and $x,y\in\bC^p$, we have
\[
<\rho_n(\bt) x,\rho_n(\bt) y>\,=\,<x,y>\,.
\]

This feature is analogous to the form of Squier \cite{Squier}
for the Burau representation $\psi_n$.

The \em{(reduced) Burau representation} $\psi_n$ of $B_n$, depending
on a complex parameter $q$ and acting on $\bC^{n-1}$, is given by:
\begin{eqnarray*}
& \psi_n(\sg_i)\,=\left[
\begin{array}{*9c}
1 \\
& \ddots & & & & & 0 \\
& & 1 \\
& & & 1 & -q \\
& & &   & -q \\
& & &   & -1 & 1 \\
& & &   & & & 1 \\
& & 0 & & & & & \ddots \\
& & &   & & & & & 1 \\
\end{array}\,\right]\quad\mbox{for $1<i<n-1$,} & \\[3mm]
& \psi_n(\sg_1)\,=\left[
\begin{array}{*5c}
-q & 0 \\
-1 & 1 & & \multicolumn{2}{c}{0} \\
& & 1 \\
\multicolumn{2}{c}{0} & & \ddots \\
& & & & 1 \\
\end{array}\,\right],\quad\mbox{and}\quad
\psi_n(\sg_{n-1})\,=\left[
\begin{array}{*5c}
1 & & & \multicolumn{1}{c}{0} \\
& \ddots \\
& & 1 \\
\multicolumn{2}{c}{0} & & 1 & -q \\
& & & 0 & -q \\
\end{array}\,\right]\,,
\end{eqnarray*}
where at position $(i,i)$ there is always the entry $-q$.

We used Squier's form
previously to carry out a study of the image of $\psi_n$ in
\cite{reiko}. Again Budney's form is definite for proper $q$ and
$t$. We became aware of Budney's result only recently, and so
we tried to adapt details we had worked out for $\psi_n$.

At least one serious obstacle was visible in advance. Due to the
quadratic increase in dimension, an argument via a rank estimate
for a simple Lie group is more complicated. Still most simple Lie
group representations have a dimension larger than quadratic (in the
rank), and with a certain amount of extra effort we will be able to
deal with them.

\begin{rem}\label{z2}
Let us here remark that when $<\,.\,,\,.\,>$ is definite, the
subgroup of linear transformations of $\bC^n$ respecting the
form is conjugate in $SL(n,\bC)$ to $U(n)$. Conversely, each such
subgroup determines the respected pairing up to complex conjugation
and multiples. This follows from that facts that the only outer
automorphism of $SU(n)$ is (entry-wise) complex conjugation, and
the centralizer of $SU(n)$ in $SL(n,\bC)$ are the scalar matrices.
\end{rem}

It will not be necessary to study the form very explicitly here.
We need only the following consequence of the formula in the
proof of Theorem 3 of \cite{Budney}. Below for a complex number
$z=re^{i\th}$ with $r\ge 0$ and $\th\in (-\pi,\pi]$
we set $|z|=r$ and $\arg z=\th$ (when $r\ne 0$).

\begin{prop}\label{xy}
If $q_i,-t_i\to 1$ with $|q_i|=|t_i|=1$ are chosen so that
$\ds\frac{1-q_i}{1-q_i\sqrt{-t_i}}\to 0$, then the Budney
form is definite at $q_i,t_i$ for $i$ large.
\end{prop}

\begin{rem}\label{z1}
The additional condition means that $\arg q_i/\arg -t_i\to 0$.
In other words, when $t_i$ is close to $-1$, one should choose
$q_i$ close to $1$ in a way depending on $t_i$.
It is clear that one can choose such $q_i,t_i$ which are
algebraically independent. This property will be used
in the proof of theorem \ref{t2}, but it will \em{not}
be relevant before \S\reference{S7}.
\end{rem}

The reason why we are interested in the value $t=-1$ is (see
\cite{Krammer}):

\begin{lemma}\label{lsym}
$\rho_n$ turns into the symmetric square of $\psi_n$ for $t=-1$.
\end{lemma}

Then, for $q=1$ we have (the symmetric square of) the permutation
homomorphism $\pi_n$. The homomorphism $\pi_n$ means here the
$n-1$-dimensional (irreducible) representation obtained from the
action of the symmetric group $S_n$ (onto which there is an obvious
homomorphism from $B_n$) on the coordinates of $\bC^n$, after
removing the (invariant) space generated by $(1,1,\dots,1)$.
The notion of symmetric square is explained more precisely below.

\begin{warn}\label{w1}
The following should be kept in mind regarding the $t$ variable.
\begin{enumerate}
\item
The convention of \cite{Krammer} for the matrices of the representation 
differs from \cite{Bigelow,Budney}; $t$ of former is $-t$ of latter.
We stick with Bigelow-Budney's convention for $t$.
\item
Also, in theorem 4.1, p.483 of \cite{Bigelow} there is a misprint: in
the fourth option $t$ should be $-t$; this is set right on p.782 of
\cite{Budney}, and our \eqref{rf} is a reproduction of latter formula.
\item
Next, $t$ is used often, e.g. in \cite{Jones}, as the variable for
$\psi_n$, but it is here $q$, \em{not} $t$, that via lemma \ref{lsym}
originates from Burau (and what we may call the Burau variable).
Apart from replacing $t$ by $q$, our definition of $\psi_n$ is as
in \cite{Jones}.
\end{enumerate}
\end{warn}

Let us clarify and fix some language. For a (complex)
vector space $V$ with basis $e_1,\dots,e_n$, we can define
the \em{symmetric square} $\Sym^2 V$ to be subspace
of $V\otimes V$ spanned by elements 
\[
v\odot w=\frac{v\otimes w+w\otimes v}{2}\,.
\]
$\Sym^2 V$ has the standard basis
\[
\{\,e_i\odot e_j\ :\ i\le j\,\}\,.
\]
When an endomorphism $f$ acts on a vector space $V$, then
it induces an endomorphism on $\Sym^2 V$ we write $\Sym^2 f$.
So we can talk of $\Sym^2 \psi_n$.

We will also need the \em{antisymmetric square} $\Lm^2 V$
generated by
\[
v\wedge w=\frac{v\otimes w-w\otimes v}{2}\,.
\]
Similarly there is a meaning to $\bigwedge^2 f$.


\begin{lemma}\label{t1_}
Assume $q,t$ with $|t|=|q|=1$ are chosen so that $t\ne \pm 1$
and $q$ is, dependingly on $t$, sufficiently close to $1$.
Also, assume the Budney form is definite at $q,t$. Then $\rho_n$
is irreducible at $q,t$.
\end{lemma}

\proof 
For $(q,t)=(1,-1)$ the symmetric square $\Sym^2\pi_n$ of the
permutation homomorphism acts by permuting the indices in $v_{i,j}$.
Now, when $q=1$, but $t\ne \pm 1$, then the action is similar,
except that $\sg_i$ acts on $v_{i,i+1}$ by multiplying by $-t$.
It is clear then that such an endomorphism has eigenvalues
$\pm 1$ and $-t$, with the eigenspace for $-t$ being
1-dimensional, generated by $v_{i,i+1}$. 

Since by unitarity every invariant subspace has an invariant
(orthogonal) complement, it follows that if there is an invariant
subspace $V$, then $V^\vd$ is also invariant, and one of both,
say w.l.o.g. $V$, contains $v_{1,2}$. But then by the way the
action is described, $V=R$.

Thus $\rho_n$ is irreducible for $q=1$, and then for $q$
close to $1$ because irreducibility is an open condition.
\qed

\begin{rem}\label{z4}
The reason we propose this proof, apart from its simplicity,
is to outline a way in which the algebraic independence
condition on $q,t$ can be circumvented. This condition inavoidably
enters if one likes to return from formal algebra to complex-valued
$q,t$, as in the approaches of \cite{Marin} and \cite{Zinno},
where in latter paper, irreducibility follows from the
identification of $\rho_n$ to a summand of the BWM algebra. 
Algebraic independence is `generically' satisfied, but for many
concrete values of $q,t$ it may be false, or at least difficult
to establish. (It is needed for the faithfulness, but
latter is not essential in our arguments, until \S\ref{S7}.)
It is clear from our proofs that the condition `$q$ close
to $1$' can, in both lemma \reference{xy} and \reference{t1_},
be made precise by a slightly more technical calculation.
The clarification for which parameters \em{exactly} $\rho_n$ is
irreducible requires far more effort, and is only subject of
ongoing work. The sole
written reference I received, only \em{a posteriori}, from
I.~Marin, is a very recent Ph D thesis of C.~Levaillant
\cite{Levaillant}. (Allegedly Bigelow has an own, unpublished,
proof.)
\end{rem}

\section{Lie groups}

\subsection{Correspondence between compact
and complex Lie groups\label{oo}}

We start by reviewing a few basic facts from Lie group theory,
which mostly occurred in our treatment of Burau \cite{reiko}.

Let $G$ be a connected compact Lie group with Lie algebra $\fg$.
Compactness implies in particular that $G$ is real,
finite-dimensional and linear reductive.
Linear reductive means for a closed subgroup $G\subset GL(n,\bC)$
that the number of connected components is finite and $G$ is
closed under conjugated matrix transposition $M\mapsto \ol{M^T}$.

A \em{linear representation} of $G$ is understood as a pair $\rho=
(V,\pi)$ made of a vector space $V$ and a homomorphism $\pi\,:\,
G\,\to\,\Aut(V)$. We will often omit $\pi$ and identify $\rho$
with $V$ for simplicity, if unambiguous. A representation is
\em{irreducible} (and will be often called \em{irrep} below)
if it has no non-trivial (i.e.\ proper and non-zero)
invariant subspaces. Linear reductiveness of $G$ implies that 
each invariant subspace of a linear representation of $G$ has a
complementary invariant subspace, so that each representation of
$G$ is \em{completely reducible} as direct sum of irreducible
representations.

To $G$ there exists a uniquely determined complex connected linear
reductive Lie group $G_\bC$, with \\[1.5mm]
\mylabel $\fg_\bC=\fg\otimes \bC$ is the Lie algebra of $G_\bC$, and\\
\mylabel $G\subset G_\bC$ as a closed subgroup. \\[1.5mm]
\setcounter{type}{0}%
Then $G_\bC$ is called a \em{complexification} of $G$.
If $G$ is simply connected, so is $G_\bC$, and then any other
connected complex Lie group with Lie algebra $\fg_\bC$ is a covering
of $G_\bC$. The real group $G$ is always a maximal compact subgroup
of $G_\bC$; we call it the \em{compact real form} of $G_\bC$.

Thus we have a one-to-one correspondence between a compact connected
real (simply-connected) Lie group and a (simply-connected) connected
linear reductive complex Lie group.
Under this correspondence to $G=SU(n)$ we have $G_\bC=SL(n,\bC)$.
The groups are connected and simply-connected (for $n\ge 2$).

The correspondence behaves well w.r.t.\ many properties. The
real form $G$ is simple, if and only if $G_\bC$ is too. (In
particular, if $G$ is semisimple, so is $G_\bC$.) For every
complex representation $\rho=(V,\pi)$ of $G$ (`complex' means
that $V$ is a complex vector space) we have an `extension'
to a representation $\tl \rho=(V,\pi_\bC)$ of $G_\bC$, such that
$\pi_\bC$ is an extension of $\pi$ from $G$ to $G_\bC$. This
extension is still faithful if $G$ is compact.

\subsection{Symmetric pairs\label{SSP}}

Let $G$ be a Lie group and $\sg$ an involution.
Define 
\[
G^\sg\,:=\{\,g\in G\,:\,\sg(g)\,=\,g\,\}\,
\]
to be the $\sg$-invariant subgroup of $G$ and $G^\sg_0$
the connected component of the identity. Then a pair
$(G,H)$ for a closed subgroup $H$ with $G^\sg_0\subset H\subset
G^\sg$ is called a \em{symmetric pair}.

In the case $G=SU(n)$ the symmetric pairs have been classified
by Cartan. See \cite[Chapter IX.4.A, table p.354]{Helgason}.
In this case $H$ is some of $S(U(m)\times U(n-m))$, $Sp(n/2)$
if $n$ is even, or $SO(n)$. 

Let us give the corresponding involutions $\sg$ that define the
symmetric pairs (see p.~348 of \cite{Helgason}). 

Define $M_{i,j}$ to be the matrix with all entries $0$ except that
at the $(i,j)$-position, which is $1$. Let $\dig(x_1,\dots,x_n)=
\sum_{i=1}^nx_iM_{i,i}$ be the diagonal matrix with entries
$x_1,\dots,x_n$, so that $Id_n=\dig(1,\dots,1)$ (with $n$ entries
`$1$') is the identity matrix. 

For $S(U(m)\times U(l))$, with $m+l=n$, the involution $\sg$
is of the form $\sg_{m,l}\,:\,M\mapsto I_{m,l}MI_{m,l}$, where 
\[
I_{m,l}\,=\,\dig(1,\dots,\inx{1}{m},\inx{-1}{m+1},\dots,-1)\,.
\]
For $n=2n'$ even, $Sp(n')$ respects the involution $\sg_J
\,:\,M\mapsto J^{-1}\bar MJ$, where
\begin{eqn}\label{Jsym}
J\,=\,\left[\begin{array}{c|c}
\rule[-4mm]{0mm}{2.9em}\rx{1mm}0\rx{1mm} & \rx{1mm}-Id_{n'}\rx{1mm}\\
 \hline \rule[-4mm]{0mm}{2.9em}\rx{1mm}\phantom{-}Id_{n'}
\rx{1mm} & \rx{1mm}0\rx{1mm}
\end{array}
\right]\,,
\end{eqn}
and $\bar M$ is the complex conjugation (of all entries) of $M$.
For $SO(n)$, the involution $\sg$ is given by $\sg(M)=\bar M$.

These subgroups can be also defined in the standard representation
by the linear transformations that respect a certain (complex)
non-degenerate bilinear form, which is Hermitian, skew-symmetric
or symmetric resp.
All transformations that respect such a form determine, up to
conjugacy, a subgroup of one of the three types.

Analogous three types of subgroups $R(m,n,\bC)$, $Sp(n/2,\bC)$ and
$SO(n,\bC)$ can be defined for $SL(n,\bC)$. Here $R(m,n,\bC)$ is
the group of all (complex-)linear unit determinant transformations
of $\bC^n$ that leave invariant a subspace of dimension $m$.
(In contrast to the unitary case, there is not necessarily a
complementary invariant subspace!)

We call the three types of groups \em{reducible}, \em{symplectic}
and \em{orthogonal} resp.
We call a representation $V$ of $G$ after one of the types, if
it is contained in a conjugate of a group of the same name.
Orthogonal and symplectic subgroups/representations
will be called also \em{symmetric}, the others \em{asymmetric}.

In the real-complex correspondence we explained, we have that
\em{if the representation $\rho$ of $G$ is symmetric, then so is the
representation $\tl \rho$ of $G_\bC$}. This is easily seen by
restricting the respected bilinear form to the reals.

\section{The maximal subgroups\label{mz}}

\begin{conv}\label{conv1}
{}From here, until the conclusion of the proof of theorem
\reference{t1} at the end of \S\reference{S5}, we will
assume, unless we clearly indicate otherwise, that
$q,t$ are fixed unit norm complex numbers that \em{satisfy
the assumptions of theorem \reference{t1}}. In this
situation, we will usually omit explicit reference to
the parameters $q,t$ in the notation.
We may mention here that these conditions are 
stronger than we actually need, and were chosen so
as to keep the formulation of theorem \reference{t1} simpler.
We will elaborate on weaker (but slightly more technical)
sufficient conditions in \S\reference{S6}.
\end{conv}


For our approach to theorem \reference{t1}, we will mainly 
study the \em{normalization} $\rho_n'\subset SU(p)$. That
is, we consider the (determinant) factorization $U(p)=SU(p)
\times U(1)$ and the projection on the first factor. This
means that we multiply $\rho_n$ by a power of a scalar
$\mu\in\bC$ (depending on $q,t$),
\[
\rho'(\bt)\,=\,\mu^{[\bt]}\,\cdot\,\rho_n(\bt)\,,
\]
so that $\rho'(\bt)$ has determinant $1$.
Here $[\bt]$ is the exponent sum of $\bt$, its image under
the homomorphism $B_n\to \bZ$ sending all $\sg_i$ to $1$.
This scalar $\mu$ can be calculated as in \eqref{qmu}.

We will prove that 
\begin{eqn}\label{d1}
\ol{\rho_n'(B_n)}=SU(p)\,.
\end{eqn}
If then $\ol{\rho_n(B_n)}\ne U(p)$, it must have codimension
one, and so be a collection of components isomorphic to
$SU(p)$. But we have
\begin{eqn}\label{dt}
\dt( \rho_n(\sg_i) )\,=\,-t(-q)^n\,,
\end{eqn}
which is not a root of unity by assumption, so the
projection of $\ol{\rho_n(B_n)}$ onto $U(1)$ is not discrete.

To prove \eqref{d1} we argue indirectly, and assume the contrary.

Since $\hH=\ol{\rho_n'(B_n)}$ is a compact Lie subgroup of $SU(p)$,
we can complexify its connected component $\hH_0$ of the identity,
and obtain a (faithful) representation of a
reductive complex Lie subgroup $\tl H$ of $SL(p,\bC)$, which is a
proper subgroup by dimension reasons. It is contained in a
maximal proper subgroup which we call $H$. We will
show that $H=SL(p,\bC)$, and have a contradiction.

\begin{conv}
Let us stipulate that until the end of the proof
of lemma \reference{pp6}, we use $n$ to indicate the dimension
rather than the number of braid strands ($\rho_n$ will not
appear in this scope).
\end{conv}

In the 1950s, Dynkin published a series of seminal papers
in the theory of Lie groups. One of his remarkable achievements
was the classification of maximal subgroups of classical Lie groups
\cite{Dynkin}. (See theorems 1.3, 1.5 and 2.1 in \cite{Dynkin}.)

\begin{theo}\label{prop0}(Dynkin \cite{Dynkin})
A maximal proper subgroup of $SL(n,\bC)$ is conjugate in $SL(n,\bC)$ to
\\
\mylabel some symmetric representation, i.e., $SO(n,\bC)$ or $Sp(n/2,
         \bC)$ when $n$ is even, \em{or} \\ \label{1.}%
\mylabel to $SL(m,\bC)\otimes SL(m',\bC)$ with $mm'=n$ and $m,m'\ge 2$
         (one which is non-simple irreducible), \em{or} \\
\mylabel to $R(m,n,\bC)$ (one which is reducible), \em{or} \\
         \label{3.}%
\mylabel it is an irreducible representation of a simple Lie group. 
	 \label{4.}%
\setcounter{type}{0}%
\end{theo}

%
%

For a non-simple group $H=H_1\times H_2$, one considers $\bC^n=
\bC^{mm'}\simeq \bC^{m}\otimes \bC^{m'}$ as a tensor (Kronecker)
product, and $H_1$ resp. $H_2$
acts on $\bC^{m}$ resp. $\bC^{m'}$. 

We call the first 3 types of subgroups \em{orthogonal},
\em{symplectic}, \em{product} and \em{reducible} resp.
We will exclude these types in the next section, before
we get to deal with \reference{4.}.


\section{Excluding symmetric, product and reducible
  representations\label{Sym}}

Note that case \reference{1.} in theorem \reference{prop0} is in fact
included in case \reference{4.}. It is singled out because it can
be handled by a more elementary eigenvalue analysis. What theorem
\reference{prop0} actually achieves is also a description of the
maximal subgroups of $SO(n,\bC)$ and $Sp(n/2,\bC)$, and one can
decide which of the representations of case \ref{4.} are symmetric.

It will be useful to have an elementary analysis of the eigenvalues
on the various subgroups.

\begin{prop}\label{pp5}
If $\lm\in \bC$ is an eigenvalue of an orthogonal or symplectic
matrix $M$ in $SU(n)$, then so is $\bar\lm$.
\end{prop}

\proof Orthogonal matrices are conjugate to real ones, so their
characteristic polynomial has real coefficients. So consider
a symplectic matrix $M$. Then the operator $J$ in \eqref{Jsym}
of the involution $\sg_J$ that $M$ respects (up to conjugacy)
satisfies $J^2=-Id_n$. So $M=J^{-1}\bar MJ$ is equivalent to
$MJ=J\bar M$. Now by assumption there is a vector
$v\in\bC^{n}\sm \{0\}$ with $Mv=\lm v$. Then 
\[
MJ\bar v=J\bar M\bar v=J\ol{Mv}=J\ol{\lm v}=\bar \lm J\bar v\,,
\]
so $J\bar v\ne 0$ is an eigenvector of $M$ for eigenvalue
$\bar \lm$. \qed

\begin{lemma}\label{pp5'}
Assume that $H=H_1\otimes H_2\in SL(n_1,\bC)\otimes SL(n_2,\bC)$.
Let $\lm_i$ be the eigenvalues of $H_1$ (counting multiplicities)
and $\mu_j$ those of $H_2$. Then the eigenvalues (counting
multiplicities) of $H$ are $\lm_i\mu_j$.
\end{lemma}

\proof $H_{1,2}$ have by Jordan box decomposition lower-triangular
bases $e_i$ and $f_j$. Order the tensor basis $\{e_i\otimes f_j\}$
so that elements with smaller $i+j$ appear first. Then w.r.t. the
so ordered basis, $H=H_1\otimes H_2$ is lower-triangular.
Since in lower-triangular matrices, the eigenvalues appear on
the diagonal, the claim is clear. \qed

\begin{lemma}\label{pp6}
Assume that for $H=H_1\otimes H_2\in U(n_1)\otimes U(n_2)$ all
eigenvalues $\lm_i$ except exactly one (and single-multi\-pli\-ci\-ty)
eigenvalue $\lm_{i_0}$ satisfy $\lm_i\in\{1,a\}$ for $a\in\bC$ with
$a\ne \pm 1$. Furthermore assume that \\
\mylabel $\lm_{i_0}\ne a^2$ and that \\
\mylabel $\lm_{i_0}\ne 1/a$, or $\lm_i=1$ for at least half of all
$\lm_i$.\\
\setcounter{type}{0}%
Then one of $n_1$ or $n_2$ is equal to $1$,
in other words, $H$ is a direct product only in a trivial way.
\end{lemma}

\proof Let $H_1$ have eigenvalues $\mu_j$ and $H_2$ have eigenvalues
$\nu_k$. Since $H_{1,2}$ are diagonalizable, $H$ has eigenvalues
$\lm_i=\mu_j\nu_k$, with all $\mu_j,\nu_k\ne 0\,.$

Let $\lm_{i_0}=\mu_{j_0}\nu_{k_0}$ and assume $n_1,n_2>1$,
where $n_1$ is now the number of $\mu_j$ and $n_2$ that of $\nu_k$.
By choosing $k\ne k_0$ and looking at the set $\{\,\mu_j\nu_k\,\}$
(for the fixed $k$ but varying $j$), we see that $\{\,\mu_j\,\}
\subset \{x,y\}$ and $x/y=a^{\pm 1}$. A similar conclusion applies
to $\{\,\nu_k\,\}$ using $\mu_j$ for some fixed $j\ne j_0$.

Then it is clear that $\lm_{i_0}\in\{a^2,1/a\}$, and we
assumed that the former value is not taken. If $\lm_{i_0}=1/a$,
then (for $a\ne \pm 1$) it follows that exactly one of the
$\mu_j$ is different from all the others, which are equal, and
similarly for $\nu_k$. This implies that the multiplicity
of $\lm_i=a$ is $(n_1-1)(n_2-1)\ge n/2-1$, for the number
$n=n_1n_2$ of eigenvalues $\lm_i$ (with $n_1,n_2>1$). Then,
$\lm_i=1$ occurs $n-(n_1-1)(n_2-1)-1\le n/2$ times, which we
excluded. This gives a contradiction to the assumption
$n_1,n_2>1$. \qed

Now we apply the preceding lemmas to $\rho_n$ (with
$n$ resuming the meaning of number of braid strands).

\begin{lemma}\label{hhi}
The image of $\rho_n'$ for $n\ge 3$ is not orthogonal or
symplectic, and hence neither is $H$.
\end{lemma}

\proof Assume first that $\rho_n'(\sg_1)\in \hH_0$,
in the notation of \S\reference{mz}.
The eigenvalues of $\rho_n(\sg_1)$ can be easily determined
from the Krammer matrix, and replacing $t$ by $-t$ according to
the warning \reference{w1}. The result is
\begin{eqn}\label{ero}
\Big\{\,-tq^2\,,\ \ \underbrace{-q,\dots,-q}_{\scbox{$n-2$ times}}\,,
\ \underbrace{1,\ \ 1,\ \ \dots,1}_{\shortstack{\scbox{$(n-1)(n-2)/2$}\\
\scbox{times}}}\,\Big\}\,.
\end{eqn}
Let us fix, also for outside this proof, the following notation.
A set of $n$ copies of $k$ will be written $\{k\}^n$, and union
of such sets will be written as a product. Then the above set
can be written as $\{-tq^2\}\{-q\}^{n-2}\{1\}^{(n-1)(n-2)/2}$.

To normalize for the eigenvalues of $\rho_n'(\sg_1)$,
this set has to be multiplied by 
\begin{eqn}\label{qmu}
\mu=\det(\rho_n(\sg_1))^{-2/n(n-1)}\,,
\end{eqn}
with \eqref{dt}. 
For the chosen $q,t$, none of the resulting numbers is
real (i.e. $\pm 1$). This in particular finishes the cases
$n(n-1)/2$ odd, so $n\ge 4$. But then there is a pre-dominant
occurrence of $\mu$, and the set is not closed under conjugation.

Now let $\rho_n'(\sg_1)\not\in \hH_0$. By the assumption on
$q,t$ in theorem \ref{t1}, no two distinct eigenvalues $\lm_i$
of $\rho_n(\sg_1)$ have a quotient which is a root of unity.
Then one can choose a number $m$ large enough so that all
eigenvalues of $\rho_n(\sg_1^m)$ are as close to $1$ as desired.
(This can be seen for example by looking at the closure of the
infinite cyclic subgroup generated by $(\lm_1,\lm_2,\lm_{3})$ for
the 3 distinct $\lm_i$ within the $3$-dimensional torus $T^3$,
and arguing that this closure, which is formally a Lie subgroup
of $T^3$, cannot have a codimension.) Therefore, by compactness
$\rho'_n(\sg_1^m)\in \hH_0$. Then one argues analogously to
above with $\rho_n'(\sg_1^m)$. \qed

\begin{lemma}\label{lip}
The image of $\rho_n'$ is not contained in a Kronecker product.
\end{lemma}


\proof Consider first the situation when $\rho_n'(\sg_1)\in \hH_0$. 
We want to show that there are no (non-trivial) matrices $H_1, H_2$ with
\begin{eqn}\label{xz}
\rho_n'(\sg_1)(q,t)=H_1(q,t)\otimes H_2(q,t)\,.
\end{eqn}
%
To rule out \eqref{xz}, we may replace $\rho_n'$ by $\rho_n$,
and consider the eigenvalues of $\rho_n(\sg_1)$.
Using \eqref{ero}, and under the restrictions on $q,t$ of convention
\reference{conv1}, we can apply lemma \reference{pp6}.
It gives the desired conclusion. 

%
%

Again, for $\rho_n'(\sg_1)\not\in \hH_0$, one argues with
$\sg_1^m$ and replaces $t$ by $t^m$ (to which the
restriction of theorem \reference{t1} applies in the
same way). \qed


\begin{lemma}\label{lred}
The group $H$ acts irreducibly on $\bC^p$.
\end{lemma}

\proof Since, in the notation of \S\reference{mz}, we have
$H\supset \tl H\supset \hH_0$, it is enough to prove that
irreducibility of $\rho'_n$, or simpler of $\rho_n$, is not spoiled
when we restrict $\hH$ to $\hH_0$.

When the Budney form is definite, each $\rho_n$-matrix
diagonalizes. Thus invariant subspaces and irreducibility
will be preserved if we pass to $m$-th powers of any generating set
$\{\tau_i\}$ of $B_n$, provided there are no two distinct eigenvalues
of any $\rho_n(\tau_i)$ which differ (multiplicatively) by a
root of unity.

By the assumption on $q,t$ in theorem \ref{t1}, this condition holds
for $\tau_i=\sg_i$ from \eqref{ero}. Clearly one can choose
$m$ in the last paragraph of the proof of lemma \reference{hhi}
so that the eigenvalues of all $\rho_n(\sg_i^m)$ for
$1\le i\le n-1$ are as close to $1$ as desired,
and then all $\rho_n'(\sg_i^m)\in \hH_0$. \qed

\section{Rank estimate\label{RE}}

In the quest for what $H$ could be, we are left from the list of
theorem \ref{prop0} only with case \ref{4.}. To deal with this,
in the following it is necessary to appeal to a larger
extent to the Lie theory described, mainly in the appendix,
in \cite{Dynkin}. We will repeat a certain part, though we
would have to refer there for further details. 

\begin{conv}
References to
pages, and to equations or statements numbered `$0.\,\cdots$' are
to be understood \em{to Dynkin's paper} (in the translated version).
\end{conv}

The \em{rank} $\rk G$ of a simple Lie group $G$ is the maximal dimension
of a torus $G$ contains, or the number of nodes in its Dynkin diagram.
The latter description will be used from the next section on.
Here we have to deal with the torus.

We will recur our rank estimate for $\rho_n$ to the one for $\psi_n$
by means of the important observation in lemma \reference{lsym}.

Let $\bt$ be a \em{fixed} braid in $B_n$. A braid $\bt\in B_k$ for
$k\le n$ can be regarded also as a braid $\bt\in B_{1,k}\subset B_n$.
The following lemma tells us how to determine the eigenvalues of
$\rho_n(\bt)$.

Let for a matrix $M$, by $E=\Ev M=\{\lm_i\}$ be denoted the
eigenvalues of $M$ (counting multiplicities), and let
$\Sym^2 E\,=\,\{\,\lm_i\lm_j\,:\,i\le j\,\}$.

\begin{lemma}\label{o1}
\[
\Ev \rho_n(\bt)\,=\,\bigl(\,\Sym^2 \Ev\psi_n(\bt)\,\sm\,
\Sym^2 \Ev\psi_k(\bt)\,\bigr)\,\cup\,\Ev \rho_k(\bt)\,.
\]
For $q,t$ of definite Budney form, the eigenspaces of $\rho_n(\bt)$
of eigenvectors in $\Ev \rho_k(\bt)$ correspond to $\cE_k:=\{
\,v_{ij}\,:\,1\le i<j\le k\,\}$.
\end{lemma}

\proof We order the basis $\cE_n$ of $\bC^p$ so that $\cE_k$
occur first. It is obvious from the definition \eqref{rf} that
$\rho_n\restr{B_k}$ respects $\cE_k \subset \cE_n$. So the
matrix of $\bt\in B_k\subset B_n$ has the form
\begin{eqn}\label{kn}
\rho_{n}(\bt)\,=\,
\def\k{\kern2mm}
\left[
\begin{array}{c|c@{\ \ }c@{\ \ }c}
\ry{2.0em}\rho_{k}(\bt) &    &    0   &   \\[3mm]%
\hline%
                  &    &   \ry{1.1em}     &   \\[-2.8mm]
A      &    & B &   \\[2.3mm]
\end{array}
\right]\ ,
\end{eqn}
where $A,B$ also depend on $\bt$. Thus
\begin{eqn}\label{q1}
\Ev \rho_n(\bt)=\Ev \rho_k(\bt)\cup \Ev B(\bt)\,.
\end{eqn}

The next important observation is that by definition
the variable $t$ does not occur in $B$ for $\rho_n(\sg_i)$,
$i=1\,\dots,k-1$. Then the same is also true for
their inverses, and finally thus for $\rho_n(B_k)$.
But since $B$ does not depend on $t$, its eigenvalues can be
determined setting $t=-1$. Then $\rho_k=\Sym^2\psi_k$
and $\rho_n=\Sym^2\psi_n$. Thus we have
\begin{eqn}\label{q2}
\Ev B(\bt)\cup \Sym^2\Ev\psi_k(\bt)\,=\,\Sym^2\Ev\psi_n(\bt)\,.
\end{eqn}
Combining \eqref{q1} and \eqref{q2}, we have the claim. \qed

Note that \eqref{ero} also follows from this lemma.

Now we apply the lemma on the following elements in $B_n$
that were of central importance also for Burau.
\begin{eqn}\label{q3}
\bt_{n,k}\,=\,\Dl_k^2\,=\,(\sg_1\dots\sg_{k-1})^k\,\in\,B_n\,.
\end{eqn}

\begin{lemma}\label{re}
When $q,t$ are chosen as in theorem \reference{t1}, and the
Budney form is definite, we have
\[
\ol{\rho_n'(\,<\,\bt_{n,2},\dots,\bt_{n,n-1}\,>\,)}\,=\,T_{n-2}\,,
\]
an $n-2$-dimensional torus. Thus in particular $\rk H\ge n-2$.
\end{lemma}

\proof We have that $\rho_k(\bt_{n,k})$ are scalars, and
with the notation explained below \eqref{ero}, we have
\[
\Ev \psi_n(\bt_k)=\{q^k\}^{k-1}\{1\}^{n-k}\,,
\]
as observed
in \cite{reiko}. If the Budney form is definite on $B_n$,
then so it is on $B_k$, and all matrices are diagonalizable.
The eigenspaces corresponding to the eigenvalues $\{q^k\}^{k-1}$
of $\psi_n$ are spun by $e_1\,\dots\,e_{k-1}$. The claim follows by
a careful look at eigenvalues and eigenspaces. \qed

\section{The irreps of simple Lie groups\label{S5}}

\subsection{Dynkin diagrams. Weyl's dimension formula}

With \S\reference{RE} we are left to consider 
irreps of simple Lie groups. Moreover we know quite
exactly the dimension. It will be more convenient to
look at our rank estimate in lemma \reference{re} from the
point of view of the group, not the number of braid
strands. It says then that for a group of rank $n$, the dimension
of the representation must be $(n'+1)(n'+2)/2$ for some $n'\le n$.
Moreover, this irrep should not have an invariant form
by \S\reference{Sym}. We try next to find out how to obtain
all these irreps.

By the work of Cartan, irreps $\phi$ of a simple Lie group $G$ are
determined by their \em{highest weight} $\Lm$, and latter
is completely described by the property that
\begin{eqn}\label{a_i}
a_i\,=\frac{2(\Lm,\ap_i)}{(\ap_i,\ap_i)}
\end{eqn}
are non-negative integers for all simple roots $\ap_i$ of 
the Lie algebra $\fg$ of $G$ (and at least one $a_i$ is positive).
The scalar product $(\,.\,,\,.\,)$
is the one defining the Dynkin diagram: with the normalization
that (in the cases we require below) all simple roots have
length $1$, nodes of the Dynkin diagram depicting orthogonal
vectors are not connected, and connected nodes correspond
to vectors of scalar product $-\myfrac{1}{2}$.

Since the $\ap_i$ correspond to nodes in the Dynkin diagram,
our convention, as in \cite[p.329 top]{Dynkin}, will be to write
$a_i$ at the node for $\ap_i$ in the diagram (but omit zero
entries). We will refer to $a_i$ also as \em{labels} of the nodes.
If $a_i=0$ for all $i$ except exactly one, where $a_i=1$, we
call $\phi$ a \em{basic} representation; it is obviously 
associated to the simple root (or node) $\ap_i$ with $a_i=1$.

The following formula calculates the dimension $N(\phi)$ 
of the irrep $\phi$ corresponding to $\Lm$.

\begin{lemma}(Weyl's formula, Theorem 0.24)
\begin{eqn}\label{We}
N(\phi)\,=\,\prod_{\ap\in\Sg_+}\,\frac{(\Lm+g,\ap)}{(g,\ap)}\,,
\end{eqn}
where $\Sg_+$ is the set of positive roots of $\fg$
and
\begin{eqn}\label{We'}
g=\frac{1}{2}\,\sum_{\bt\in \Sg_+}\,\bt\,.
\end{eqn}
\end{lemma}

Mostly we will appeal to the following consequence.
(We have used the fact in a weaker form already for Burau.)

\begin{lemma}\label{ij}
Assume that one increases the label $a$ of a node in the Dynkin
diagram to $a+1$ ($a\ge 0$). Then the dimension of the
irrep grows a least by a factor of $(a+2)/(a+1)$, in particular
it increases strictly.
\end{lemma}

\proof If $\ap$ is a simple root, then it is known that 
\[
\frac{2(g,\ap)}{(\ap,\ap)}\,=\,1\,,
\]
(see (0.141)),
so that $(g,\ap)>0$. Then this is also true for all $\ap\in \Sg_+$,
since they are just sums of (distinct) simple roots. For the same
reason $(\Lm,\ap)$ are just sums of $a_i$, and so non-negative.

Thus increasing some $a_i$ will not decrease any of the (positive)
factors in the product of \eqref{We}. The estimate of increase of
$N(\phi)$ follows from looking at the factor that corresponds to
$\ap_i$.
\qed

\begin{defi}
Let us say that an irrep $\phi'$ \em{dominates} another
irrep $\phi$, if the labels $a_i'$ of $\phi'$ and
$a_i$ of $\phi$ satisfy $a_i'\ge a_i$ for all $i$.
\end{defi}

\begin{figure}
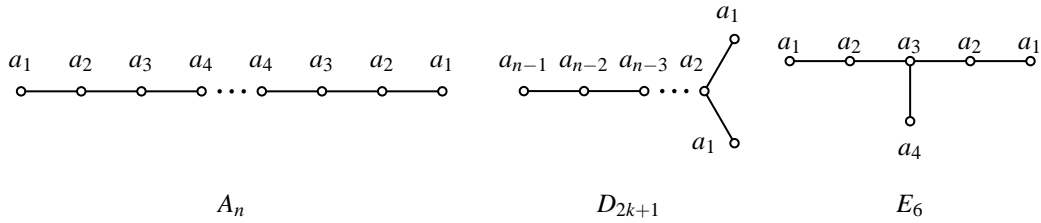

\[
\begin{array}{c@{\quad\qquad}c@{\qquad}c}
  \diag{8mm}{7}{1}{
   \pictranslate{0 0.2}{
    \picline{0 0.3}{3 0.3}
    \picline{4 0.3}{7 0.3}
    \picmultigraphics{3}{0.2 0}{\picfilledcircle{3.3 0.3}{0.02}{}}
    \picmultigraphics{4}{1 0}{\picfilledcircle{0 0.3}{0.07}{}}
    \picmultigraphics{4}{1 0}{\picfilledcircle{4 0.3}{0.07}{}}
    \picputtext{0 0.8}{$a_1$}
    \picputtext{7 0.8}{$a_1$}
    \picputtext{1 0.8}{$a_2$}
    \picputtext{6 0.8}{$a_2$}
    \picputtext{2 0.8}{$a_3$}
    \picputtext{5 0.8}{$a_3$}
    \picputtext{3 0.8}{$a_4$}
    \picputtext{4 0.8}{$a_4$}
   }
  }
  &
  \diag{8mm}{3.5}{2.5}{
   \pictranslate{0 0.55}{
    \picline{0 0.7}{2 0.7}
    \picline{3.5 0.7 3 q 2 : +}{3 0.7}
    \picline{3.5 0.7 3 q 2 : -}{3 0.7}
    \picmultigraphics{3}{0.2 0}{\picfilledcircle{2.3 0.7}{0.02}{}}
    \picmultigraphics{4}{1 0}{\picfilledcircle{0 0.7}{0.07}{}}
    \picfilledcircle{3.5 0.7 3 q 2 : -}{0.07}{}
    \picfilledcircle{3.5 0.7 3 q 2 : +}{0.07}{}
    \picputtext{3.4 1.14 3 q 2 : +}{$a_1$}
    \picputtext{3 0.7 3 q 2 : -}{$a_1$}
    \picputtext{2.8 1.2}{$a_2$}
    \picputtext{2 1.2}{$a_{n-3}$}
    \picputtext{1 1.2}{$a_{n-2}$}
    \picputtext{0 1.2}{$a_{n-1}$}
   }
  }
  &
  \diag{8mm}{4}{2}{
   \pictranslate{0 -0.2}{
    \picline{0 1.7}{4 1.7}
    \picline{2 1.7}{2 0.7}
    \picmultigraphics{5}{1 0}{\picfilledcircle{0 1.7}{0.07}{}}
    \picfilledcircle{2 0.7}{0.07}{}
    \picputtext{0 2}{$a_1$}
    \picputtext{4 2}{$a_1$}
    \picputtext{1 2}{$a_2$}
    \picputtext{3 2}{$a_2$}
    \picputtext{2 2}{$a_3$}
    \picputtext{2 0.23}{$a_4$}
   }
  } \\[1.1cm]
  A_n & D_{2k+1} & E_6
\end{array}
\]
\caption{The marking of symmetric irreps (i.e. which admit
an invariant symmetric or antisymmetric form).\label{FSym}}
\end{figure}

A further important tool is the decision which irreps are
asymmetric. This goes back to work of Malcev, and
can be done as explained in Theorem 0.20 p.336
and Remark C.a. on p.254. With the exclusion of symmetry in
\S\reference{Sym}, we are left only with representations of $A_n$,
$D_{2k+1}$ and $E_6$, whose labelings do \em{not} admit a certain
symmetry as shown on Figure \reference{FSym} (which reproduces
Table 3, p.365 in \cite{Dynkin}).


\subsection{$A_{n}$}

\def\Alabel#1#2#3{
  \diag{8mm}{#1}{1}{
    \picline{0 0.3}{#2 0.3}
    \picline{#2 1 + 0.3}{#2 1.5 + 0.3}
    \picmultigraphics{3}{0.2 0}{\picfilledcircle{#2 .3 + 0.3}{0.02}{}}
    \picmultigraphics{#2}{1 0}{\picfilledcircle{0 0.3}{0.07}{}}
    \picfilledcircle{#2 1.5 + 0.3}{0.07}{}
    #3
  }
}
\def\Alabelone#1#2#3#4{\Alabel{#1}{#2}{
  \picputtext{#3 0.7}{$#4$}}
}
\def\Alabeltwo#1#2#3#4#5#6{\Alabel{#1}{#2}{
  \picputtext{#3 0.7}{$#4$}
  \picputtext{#5 0.7}{$#6$}}
}
\def\Afivelabelone#1#2{
  \diag{8mm}{4}{1}{
    \picline{0 0.3}{4 0.3}
    \picmultigraphics{5}{1 0}{\picfilledcircle{0 0.3}{0.07}{}}
    \picputtext{#1 0.7}{$#2$}
  }
}

\def\Afivelabeltwo#1#2#3#4{
  \diag{8mm}{4}{1}{
    \picline{0 0.3}{4 0.3}
    \picmultigraphics{5}{1 0}{\picfilledcircle{0 0.3}{0.07}{}}
    \picputtext{#1 0.7}{$#2$}
    \picputtext{#3 0.7}{$#4$}
  }
}

Let $\pi_1$ be the elementary representation
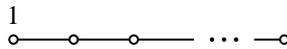
\[
\Alabelone{4.5}{3}{0}{1}
\]
of $A_n$ as $SL(n+1,\bC)$. Our assumption was that
$H\ne \pi_1$, so we will discuss the other possibilities.

First we look at the dimension of the basic representations in Table
30, p.378. These representations were written as $\pi_k$ in (0.92).
For symmetry reasons it makes sense to consider only $n\ge 2k-1$.

If $k\ge 3$, then the only case of
\[
N(\pi_k)\,=\,{n+1\choose k}\,\le\,\frac{(n+1)(n+2)}{2}
\]
is that of $k=3, n=5$. This representation
\[
\Afivelabelone{2}{1}
\]
is symmetric, and it has improper dimension 20. So we must consider
(to avoid symmetry) representations dominating
\begin{eqn}\label{111}
\Afivelabeltwo{2}{1}{0}{1}
\qquad\mbox{and}\qquad
\Afivelabeltwo{2}{1}{1}{1}
\end{eqn}
But latter two dimensions are by lemma \reference{ij}
at least $2\cdot 20=40>21$. (The exact dimensions are
$105$ and $210$, resp.)

The representations
\[
\Alabelone{4.5}{3}{1}{2}
\qquad\mbox{and}\qquad
\Alabeltwo{4.5}{3}{1}{1}{0}{1}
\]
have dimensions 
\[
\frac{(n+2)(n+1)^2n}{12}\qquad\mbox{and}
\qquad\frac{(n+2)(n+1)n}{3}\,.
\]
They are too large for $n\ge 3$. The dimension of
\[
\Alabelone{4.5}{3}{0}{2}
\]
is precisely $\frac{(n+1)(n+2)}{2}$, so that the first node
label cannot be increased beyond 2.

There remain the irreps
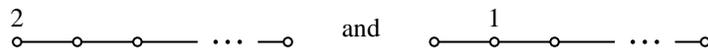
\[
\Alabelone{4.5}{3}{0}{2}
\qquad\mbox{and}\qquad
\Alabelone{4.5}{3}{1}{1}
\]
These are $\Sym^2\pi_1$ and $\bigwedge^2\pi_1$.
This can be seen as follows: the inclusions are present
by p.342 l.-8 in the first case, and p.347 l.11
in the second case. Equality follows by calculating dimensions
of the irreps using (0.148).

So it remains to rule out the symmetric and antisymmetric
square of the elementary representation $\pi_1$.
For symmetry reasons $\Sym^2\pi_1$ must be considered for
$n\ge 2$, and $\bigwedge^2\pi_1$ for $n\ge 4$ only.
By comparing dimensions, it is clear that 
the number $m$ of strands of $\rho_m$ is $m=n+2\ge 4$
in the first, and $m=n+1\ge 5$ in the second case.

We will have to
count again eigenvalues, and use the following suggestive

\begin{lemma}\label{yu}
If $M$ has eigenvalues $\lm_i$ (counting multiplicities)
then the eigenvalues of $\Sym^2 M$ are $\{\,\lm_i\lm_j\,:\,i
\le j\,\}$ and those of $\bigwedge ^2 M$ are $\{\,\lm_i\lm_j
\,:\,i< j\,\}$.
\end{lemma}

\proof Similar to lemma \reference{pp5'} (only this time with
bases $e_i\odot e_j$ resp. $e_i\wedge e_j$). \qed

We can thus finish the $A_n$ series with the following 

\begin{lemma}
$\rho_n'(\sg_1^m)$ is not a symmetric square for $n\ge 4$
and not an antisymmetric square for $n\ge 5$.
\end{lemma}

\proof It is clear that one can argue with $\rho_n$ instead.
Consider again $q=-1$. Then all eigenvalues but exactly one are $1$.
Let $\lm_i$ be the ($n-1$ resp. $n$) eigenvalues of a potential
matrix, whose (anti)symmetric square is $\rho_n(\sg_1)$.

For symmetric square $n-1$ of the $p=n(n-1)/2$ eigenvalues
must be $\lm_i^2$, and the rest $\lm_i\lm_j$ for $i<j$.
So $\lm_i^2=1$ for all but at most one $i$.
If at least three $\lm_i=\pm 1$ occur, different signs
are easily ruled out. The option $n=4$, $\lm_1=1$,
$\lm_2=-1$ and $\lm_3\ne \pm 1$ is also easy to exclude.
Clearly not all $\lm_i=\pm 1$, so all but
exactly one are (with the same sign). But then again
we must have all $\lm_i\lm_j=1$ for $i\ne j$, which is impossible.

For antisymmetric square the eigenvalues $\lm_i\lm_j=1$ for
$i\ne j$ are $1$ except one. Let us permute indices $j$ of
$\lm_j$ so that the exceptional one $-t$ occurs as $\lm_1\lm_2$.
Then all $\lm_3\lm_j$ are equal, so either all $\lm_j$
are equal, or equal except one ($\lm_3$). Neither option is
possible. \qed

\subsection{$D_{2k+1}$}

Let $n=2k+1\ge 5$ be the rank.
{}From Figure \reference{FSym} it is evident that
in order the irrep to lack symmetry, we need to label
non-trivially some of the extreme nodes labeled in figure
\reference{FSym} as $a_1$.
They correspond to basic representations called \em{spinor
representations} (see p.351). According to Table 30 page 378, their
dimension is $2^{n-1}$. We call the nodes $a_1$ thus below
\em{spinor nodes}.

Now $2^{n-1}\le (n+1)(n+2)/2$ only if $n=5$.
In this case the only nodes with associated basic
representations of dimension at most $(n+1)(n+2)/2=21$
are the spinor node, and the other extreme node. So we are
left with considering representations dominating
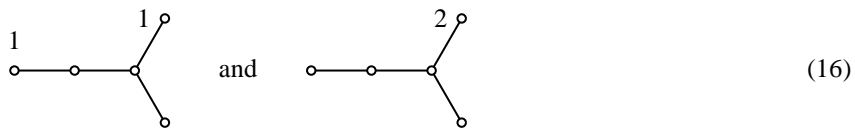
\begin{eqn}\label{D5}
\def\Dfivelabel#1#2#3{
  \diag{8mm}{2.5}{2}{
   \pictranslate{0 0.3}{
    \picline{0 0.7}{2 0.7}
    \picline{2.5 0.7 3 q 2 : +}{2 0.7}
    \picline{2.5 0.7 3 q 2 : -}{2 0.7}
    \picmultigraphics{3}{1 0}{\picfilledcircle{0 0.7}{0.07}{}}
    \picfilledcircle{2.5 0.7 3 q 2 : -}{0.07}{}
    \picfilledcircle{2.5 0.7 3 q 2 : +}{0.07}{}
    \picputtext{2.15 0.7 3 q 2 : +}{$#1$}
    \picputtext{#2 1.2}{$#3$}
   }
  }
}
\Dfivelabel{1}{0}{1}
\qquad\mbox{and}\qquad
\Dfivelabel{2}{0}{\ }
\end{eqn}
But by the estimate in lemma \reference{ij} their dimensions
are at least $3\cdot 16/2=24$. The exact computation using the
formula (0.150) gives $N(\phi)=144$ and $126$ resp.

\subsection{$E_6$}

In case of $H$ being an $E_6$ representation, the possible
dimensions are $(n+1)(n+2)/2$ for $n\le 6$. A look at the
dimensions of the basic representations of $E_6$ in Table 30,
p.378 shows that only the one marked 27 still fits this bound,
and since dimensions do not match, we must try the irrep
\begin{eqn}\label{E62}
  \diag{8mm}{4}{2}{
    \picline{0 1.5}{4 1.5}
    \picline{2 1.5}{2 0.5}
    \picmultigraphics{5}{1 0}{\picfilledcircle{0 1.5}{0.07}{}}
    \picfilledcircle{2 0.5}{0.07}{}
    \picputtext{0 1.0}{$2$}
  }\es.
\end{eqn}
It (and dominating representations) can be excluded from the
estimate in lemma \reference{ij}.

Originally we were aware, though, just of the strict increase
property in the lemma.
So we tried to compute the dimension of the irreps of $E_6$ using the
formula (0.153). But this revealed that the preceding calculation of
the numbers involved in the formula lacks several explanations and
has many errors. So here we provide a correction.

We consider an irrep of $E_6$ highest weight $\Lm$ with $a_i$
in \eqref{a_i} corresponding to nodes of the Dynkin diagram thus:
\[
  \diag{8mm}{4}{1.8}{
   \pictranslate{0 -0.4}{
    \picline{0 1.7}{4 1.7}
    \picline{2 1.7}{2 0.7}
    \picmultigraphics{5}{1 0}{\picfilledcircle{0 1.7}{0.07}{}}
    \picfilledcircle{2 0.7}{0.07}{}
    \picputtext{0 2}{$a_1$}
    \picputtext{4 2}{$a_5$}
    \picputtext{1 2}{$a_2$}
    \picputtext{3 2}{$a_4$}
    \picputtext{2 2}{$a_3$}
    \picputtext{1.7 0.53}{$a_6$}
   }
  }\ .
\]
First, on p.354, l.-11, $\lm_5$ should be $\lm_6$. On p.355, l.8,
one should refer to (0.95) and (0.96) instead of (0.108).
On (0.138') of page 355, way may clarify that if we normalize
the scalar product so that $||\ap_i||=1$, then $K=\myfrac{1}{12}$.
Then, on p.378, Table 29, the left node $\ap_6$ should be labeled
as $\ap_6=\lm_4+\lm_5+\lm_6+\lm$. The positive roots
$\Sg_+(E_6)$ can be then obtained from the list of roots
in (0.133) by choosing therein $q>p$ (in the first shape) and
the $+$ in either $\pm$ signs (for the second and third shapes).

The formulas on p.358, l.8-10, are almost entirely wrong. The
quantities $g_k$ and $g$ on the right are not properly explained,
but it is suggestive that the decomposition, similar to (0.145'),
\[
g=\sum_{i=1}^6 g_i\lm_i+g_0\lm
\]
for the element $g$ of \eqref{We'} or (0.140) is meant, where we
replaced the
$g$ of p.358, l.8-10 by $g_0$ to avoid confusion. (So $g$ is a vector
for us, given in (0.140), and $g_0$ is a scalar.) We have to assume
similarly to (0.146') that
\[
\sum_{i=1}^6 g_i=0\,,
\]
in order to have formula (0.153) working properly, and we should
(with our convention) replace $\myfrac{1}{2}\,g$ by $\myfrac{1}{2}
\,g_0$ therein. The formulas on p.358, l.8-10 should read then
\[
\begin{array}{l@{\qquad}l}
\ds l_k\,=\,l_6+\sum_{i=k}^5a_i, & \ds g_k=\frac{7-2k}{2}
  \mbox{\quad for $k\le 5$;}
\\[5mm]
\ds l_6\,=\,-\frac{1}{6}\,\sum_{i=1}^5\,i\,a_i\,, &
  \ds g_6=-\frac{5}{2}\,;
\\[5mm]\pagebreak
l=a_1+2a_2+3a_3+2a_4+a_5+2a_6\,, & g_0=11\,.
\end{array}
\]
To make the presentation more self-contained, we rewrite
the dimension formula (0.153) here. With
\[
m_k=l_k+g_k\,,\mbox{\quad and\quad}\,m=l+g_0\,,
\]
we have
\[
N(\phi)\,=\,\frac{m}{g}\,\prod_{1\le p<q\le 6}\frac{m_p-m_q}{g_p-g_q}
\,\prod_{1\le p<q<r\le 6}\frac{m_p+m_q+m_r+\myfrac{1}{2}\,\!m}
{g_p+g_q+g_r+\myfrac{1}{2}\,\!g_0}\,.
\]
Regardless of the aforementioned errors, the data for the dimensions
of the basic representations of $E_6$ in Table 30, p.378, are correct.
We then also found the exact dimension of the irrep \eqref{E62}
to be 351. (This, and some of the preceding, exact computations
of dimension will become helpful below.)

With this the proof of theorem \reference{t1} is complete.

\section{Generalizations\label{S6}}

After we found our proof of theorem \reference{t1}, we became
aware of the paper by I.~Marin \cite{Marin}. It turns out that
our theorem \ref{t1} is more or less equivalent to his theorem
B and its corollary in \S5.1 for groups of type $A$. His result
implies (most of) theorem \reference{t1} as follows. In the
terminology of our proof, one can conclude from Zariski density
(instead of going through Dynkin's list) that $\tl H=SL(p,\bC)$
when $t$ and $q$ are algebraically independent. (So this
is a slightly stronger restriction than ours.) The converse
implication is also quite obvious. Marin's result applies also
to other generalizations of $\rho_n$ (of his types $D$ and $E$),
and uses an entirely different description of the representation.
The proof is quite abstract and consists in looking at the Lie
algebra of our $H$.

In contrast, our proof is more direct and gives some new insight.
For example, since the eigenvalue argument can be carried out on
$\rho_n(\sg_1^2)$ instead of $\rho_n(\sg_1)$, theorem \ref{t1}
holds by replacing $B_n$ by any subgroup which contains the elements
$\bt_{n,k}$ ($k=2,\dots,n-1$) in \eqref{q3}, and on which $\rho_n$ is
irreducible. Here is a further variation, in which we made also
some effort to extract what conditions on $t$ and $q$ are really
needed in our arguments.

\begin{theorem}\label{t11}
Fix an integer $m\ne 0$. Assume $q,t$ with $|t|=|q|=1$ are chosen
so that $t$, $q$, $tq$, $tq^2$ and $tq^n$ are not roots of unity.
Assume the Budney form is definite at $q,t$, and $G\subset B_n$
is a subgroup as specified below, such that $\rho_n$ is irreducible
on $G$ at $q,t$.
\begin{enumerate}
\item If $G$ contains $<\,\sg_{2k-1}^{2m}\,>$ ($k\le n/2$),
and provided $tq^3$ is not a root of unity when $n=4$,
then $\ol{\rho_n(G)}\simeq U(p)$ (for $p=n(n-1)/2$).
\item If for fixed $a\ge 2$ and $l$, the group $G$ contains
$<\,\sg_{ak+l}^{2m}\,>$ (for all $k$ with $1\le ak+l\le n-1$), then
for $n$ large (in a way dependent on $a$, but independent on $q,t,l$
or $G$) we have $\ol{\rho_n(G)}\simeq U(p)$.
\end{enumerate}
\end{theorem}

\proof The condition on $tq^n$ is needed to reduce the problem
from $U(p)$ to $SU(p)$. The eigenvalue argument remains the same: as
long as $t$ is not a root of unity, we can get disposed of symmetry
and Kronecker product by looking at $\rho_n(\sg_1^{2mm'})\in\hH_0$
for proper $m'$. The condition on $q$, $tq$ and $tq^2$ enters in
order to keep lemma \reference{lred} working. With these
restrictions, the condition on $tq^3$ is what remains 
from the second listed assumption in lemma \reference{pp6},
which is needed to adapt the argument for lemma \reference{lip}.
(When $n>4$, then \eqref{ero} shows that the eigenvalue $1$
occurs too often.) The need to exclude these quantities from
being roots of unity (rather than just $\pm 1$) comes again from
the possibility that $\hH$ is not connected (i.e. $\hH\ne\hH_0$).

For the second claim a torus (within $H$) of dimension a
positive multiple of the number of braid strands is found
from looking at the action of $\rho_n(\sg_{ak+p}^{2m})$ on
subspaces of $v_{ak+p-1,ak+p}$ (where the condition on $tq^2$
is needed). Such a torus keeps an irrep analysis still
manageable. $E_6$ is relevant only for finitely many
$n$, and $D_{2k+1}$ needs (in order to prevent symmetry)
a spinor node marked, with an exponential increase in
dimension. Finally for $A_n$ one remarks that any other
labeling than the ones we studied would give a dimension
of the irrep, which is a polynomial in $n$ of degree $>2$.
Thus only finitely many $n$ would be relevant. 

In the case $a=2$ of the first claim, we have for a rank-$n$-%
group an irrep of dimension $n'(n'-1)/2$ for $n'\le 2n+1$;
in particular the dimension is at most $(2n+1)n$. 
A similar but slightly more involved discussion in cases,
as for the proof of theorem \reference{t1}, shows that in fact
under this weaker condition, still no irreps occur.

To conclude this it is helpful to use the dimension formulas
given in the proof of theorem \reference{t1} (rather than just
the rough estimates). We give just a few details. 

For $E_6$ the only new irrep fitting the dimension bound is
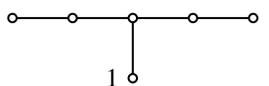
\[
  \diag{8mm}{4}{2}{
    \picline{0 1.7}{4 1.7}
    \picline{2 1.7}{2 0.7}
    \picmultigraphics{5}{1 0}{\picfilledcircle{0 1.7}{0.07}{}}
    \picfilledcircle{2 0.7}{0.07}{}
    \picputtext{1.66 0.7}{$1$}
  }\es,
\]
but it is symmetric.

For $D_{2k+1}$, after applying suitably lemma \reference{ij},
the only new possibility is the irrep of $A_7$ obtained from
the diagram on the right of \eqref{D5} by adding two nodes on the left.
But from the latter dimension calculated below \eqref{D5},
we conclude that the dimension is $>105$.

For $A_n$, the discussion is slightly lengthier. First, by using
the dimension formulas for the basic representations in Table 30,
and lemma \reference{ij}, one sees that one only needs to look
at representations where only the leftmost 3 nodes may obtain
a non-trivial label. These are discussed case-by-case.

Most of the options were studied already in the proof of theorem
\reference{t1}. We give a little information on the remaining ones,
by noticing that for the irreps
\[
\Alabeltwo{4.5}{3}{2}{1}{1}{1}
\qquad\mbox{and}\qquad
\Alabeltwo{4.5}{3}{2}{1}{0}{1}
\]
(which occur in the extended treatment of \eqref{111}),
the dimensions are $(n+1)\cdot{n+2\choose 4}$ and
$3{n+2\choose 4}$, resp. (for $n\ge 5$). The irreps 
\[
\Alabelone{4.5}{3}{0}{3}\es=\es \Sym^3\pi_1
\qquad\mbox{and}\qquad\es
\Alabelone{4.5}{3}{2}{1}\es=\es\bigwedge\nolimits^3\pi_1\,,
\]
are again (for small $n$, where the dimension estimate fails)
most conveniently ruled out by an eigenvalue argument.
\qed


The following consequence was motivated by a similar result in
\cite{Marin}. Our advantage is that our restrictions of $q,t$ are
weaker, and more explicit. (We do not appeal to the result
of Crisp-Paris either.)

\begin{corr}\label{t12}
Let $n\ge 3$.
Assume $q,t$ with $|t|=|q|=1$ are chosen so that $t$, $q$, $tq$,
$tq^2$, $tq^3$ (latter only for $n=4$), and $tq^n$ are not roots
of unity, and the Budney form is
definite and $\rho_n$ is irreducible at $q,t$. Let $m\ne 0$ be
any integer. Then 
we have 
\[
\ol{\rho_n(<\,\sg_{k}^{2m}\,:\,1\le k\le n-1\,>)}\simeq U(p)\,.
\]
\end{corr}

\proof The argument in the second paragraph of the
proof of lemma \reference{lred} with $\tau_i=\sg_i$
explains why the irreducibility of $\rho_n$ implies
the one of its restriction to the specified subgroup.
\qed

With such an argument one can treat also the \em{Hilden subgroup}
\cite{Hilden} $H_{2n}\subset B_{2n}$; from the presentation in
\cite[\S 5]{Tawn} one exhibits it to contain the elements
$\sg_{2i-1}$ ($1\le i\le n$). 
For irreducibility one needs a few extra arguments, which we provide.

\begin{prop}
Let $n>2$ be even. For $q,t$ as in lemma \reference{t1_},
$\rho_n$ is irreducible on $H_{n}$.
\end{prop}

\proof We repeat the proof of lemma \reference{t1_}, until before
the conclusion $V=R$. Now, for $q=1$, we have $R=V_1\oplus V_2$,
where $V_1$ is (linearly) generated by $v_{2i-1,2i}$, $i=1,\dots,
n/2$, and $V_2$ by all the other $v_{i,j}$. It is easily observed
that $V_k$ for $k=1,2$ are irreducible over $H_{n}$. 

If now $\rho_n(t,q)$ are reducible for $q$ converging to
$1$, then by orthogonal approximation (see \cite{reiko}) the
irrep decomposition of $\rho_n(t,q)$ is of the form
$R=V_1(t,q)\oplus V_2(t,q)$, where $V_k(t,q)\to V_k(t,1)=
V_k$ for $q\to 1$ in the sense that there are bases that
converge vector-wise; in particular $\dim V_k(t,q)=\dim V_k$.

Now again $-tq^2$ is a unique eigenvalue of
$\rho_n(t,q)(\sg_{2i-1})$ with eigenspace spun by
$v_{2i-1,2i}$. Since the matrices of $\rho_n(\sg_{2i-1})$
are conjugate, we see that some $V_k(t,q)$ must contain all
$v_{2i-1,2i}$, $i=1,\dots,n/2$, and so $V_1$. By convergence
we can have only $V_1\subset V_1(t,q)$, and by dimension
reasons $V_1=V_1(t,q)$. But it is direct to verify that for
$q\ne 1$, $V_1$ is not an invariant subspace of $\rho_n(t,q)$.
\qed

It should be remarked that \em{not} necessarily the same
$q,t$ as in lemma \reference{t1_} would do, and that the
above indirect argument spoils our control on how $q$
must be to $1$, the way we had it in remark \ref{z4}.
Still it seems not worthwhile to enter into technical
calculations in order to have this shortcoming removed,
and lemma \reference{t1_} remains at least qualitatively
true.

R.~Budney has observed irreducibility of $\rho_n\raisebox{-1mm}
{$\big|_{H_{n}}$}$ previously, at least for small $n$, but it is
the lack of written record that motivated us to supply the preceding
proposition.

\begin{corr}
Let $n\ge 4$ be even.
Assume $q,t$ with $|t|=|q|=1$ are chosen so that $t$, $q$, $tq$,
$tq^2$, $tq^3$ (if $n=4$), and $tq^n$ are not roots of unity, $t$
is close to $-1$, and $q$ is close to $1$ depending on $t$. Then
we have $\ol{\rho_n (H_{n})}\simeq U(p)$\,. \qed
\end{corr}

There seems no principal obstacle to apply our approach
to more general Artin groups, if more explicit (matrix)
descriptions of the representations are available.

\section{Non-conjugate braids\label{S7}}

The final section is devoted to the proof of Theorem \ref{t2}.

Theorem \ref{t2} was obtained first by Shinjo for knots $L$. However,
her method cannot be pleasantly applied to links, and this was our
motivation for a different approach in \cite{reiko}. We extended
Shinjo's result, showing theorem \reference{t2} when an $n-1$-braid
representation of $L$ has a non-scalar Burau matrix. One could hope
to further remove braids in the Burau kernel (which exists at least
for $n\ge 5$ \cite{Bigelow,LP}), replacing $\psi_n$ by the faithful
representation $\rho_n$. This was the origin for our interest in
$\rho_n$ in this paper. In contrast, the faithfulness of $\rho_n$
was not essential for theorem \ref{t1}. (Our approach there was
set out to apply also for many values of $q,t$ which are not
algebraically independent, and thus for which $\rho_n$ may not
be faithful.)

We should now choose some parameters $q,t$ for
which $\rho_n$ is faithful. They will have to be close to
$(1,-1)$ in the way that will get clear below, but apart from
that they should be kept fixed.

Throughout this section, $\bt\in B_{n-1}$ is a fixed non-central
braid representation of the link $L$. It will turn out very
helpful to take advantage of our work in \cite{reiko} and assume,
by having dealt with the other cases, that $\psi_{n-1}(\bt)$
is a scalar matrix. We write as
\[
C\,:=\,\{\,\ap\bt\ap^{-1}\,:\,\ap\in B_{n-1}\,\}
\]
the conjugacy class of $\bt$ in $B_{n-1}$. An element in
$C$ will typically be written as $\bt'$. Such a $\bt'$ will
be regarded also as element of $B_n$ using the inclusion
$B_{n-1}\simeq B_{1,n-1}\subset B_n$.

It is known that the center of $B_n$ is generated by the
full twist braid $\Dl^2=(\sg_1\dots \sg_{n-1})^n$.

\begin{lemma}\label{l23}
Assume for $\gm\in B_n$, that $\rho_n(\gm)$ is scalar. Then
$\gm$ is a power of the full twist braid.
\end{lemma}

\proof Scalar matrices are central, and by the faithfulness
of $\rho_n$, so must be then $\gm$. \qed

A \em{linear function} $f$ defined on the set $M(p,\bC)$ of
complex $p\times p$ matrices $M=(m_{ij})$ is an expression of the form
\begin{eqn}\label{lz}
f(M)\,=\,\sum_{i=1}^p\,\sum_{j=1}^p\,a_{ij}m_{ij}\,,
\end{eqn}
for fixed $a_{ij}\in\bC$. We call $f$ a \em{trace multiple}
if $a_{ij}=0$ and $a_{ii}=a_{jj}$ for $1\le i<j\le p$.

It is well-known that central matrices in $SU(p)$ are scalar and
that the trace is a conjugacy invariant. We showed in \cite{reiko}
that, apart from these trivial cases, there are no linear functions
of matrices invariant on a conjugacy class.

\begin{prop}\label{prop2}(\cite{reiko})
Assume that $f\,:\,M(p,\bC)\to \bC$ is a linear
function, which is not a trace multiple. Let $X$ be a non-central
element in $SU(p)$. Then $f$ is not constant on the conjugacy
class of $X$ in $SU(p)$ (considered as a subset of $M(p,\bC)$).
\end{prop}

\begin{lemma}\label{Q2}
Assume that $\psi_{n-1}(\bt)$ is a scalar matrix. Then
$\tr \rho_{n}(\bt'\sg_{n-1})$ for $\bt'\in C$ can be expressed
as a linear function of $\rho_{n-1}(\bt')$ for $-t,q$ close
to $1$. Moreover, this linear function is not a trace multiple.
\end{lemma}

\proof We assume that $q,-t$ are chosen close to $1$,
so that $\rho_{n}$ is unitary.

We note from \eqref{rf} that $\rho_{n}\restr{B_{n-1}}$ preserves
the subset $\cV_{n-1}=\{\,v_{ij}\,:\,1\le i<j\le n-1\,\}$.
By unitarity, the vectors $v_{i,n}$ for $1\le i<n$ can be
modified to $\tl v_{i,n}$, so that $\rho_{n}\restr{B_{n-1}}$
acts invariantly on the linear span $\cV'$ of $\tl v_{i,n}$.
We denote by $\tl\rho_n$ the restriction of $\rho_{n}$ (regarded
as a representation of $B_{n-1}$) to this space $\cV'$.

Since we are interested in evaluating the trace, we have the
freedom to change basis. In the basis of $\cV=\cV_{n-1}\cup
\cV'$ then we have the form \eqref{kn} (for $k=n-1$) with $B=0$.

Next we look at the matrix $A$ in \eqref{kn}. As in lemma \ref{o1}
and its proof, we noticed that the eigenvalues of $\tl\rho_{n}(\bt)$
do not depend on $t$ (although $\tl\rho_{n}(\bt)$ itself would).
So let $t=-1$. In this case we use lemma \ref{lsym} and the standard
fact (see e.g. Note 5.7 and above Example 3.2 in \cite{Jones})
that $\psi_n\restr{B_{n-1}}$ is the sum of $\psi_{n-1}$ and a
trivial representation $\tau_{n-1}$, to conclude that 
\begin{eqn}\label{zq}
\tl\rho_n=\psi_{n-1}\oplus \tau_{n-1}\,.
\end{eqn}
By multiplying the matrix of $\tl v_{i,n}$ by a proper unitary
matrix independent on $t$, we can assume that a basis of $\cV'$
is chosen w.l.o.g. so that the direct sum in \eqref{zq} is visible
in the block $A$ of \eqref{kn} for $t=-1$.

Now we assumed that $\psi_{n-1}(\bt)$ is scalar, so all its
eigenvalues are the same. But lemma \reference{o1} argued
that they do not depend on $t$, so they will be the same
also when $t\ne -1$. Since still $\tl\rho_{n}$ is unitary
for the chosen $q,t$, we see that $\tl\rho_{n}(\bt)$ is a
diagonal matrix independent on $t$ for such $q,t$, and it is
the same matrix $A=\psi_{n-1}(\bt')$ for all $\bt'\in C$.

This means that, in the basis $\cV$, the only entries of
$\rho_n(\bt')$ that vary with $\bt'\in C$ are those in
the block $\rho_{n-1}(\bt')$ in \eqref{kn}. By writing
$\rho_{n}(\sg_{n-1})$ in the same basis $\cV$, we can
then express $\tr \rho_{n}(\bt'\sg_{n-1})$ as a linear
combination of entries of $\rho_{n-1}(\bt')$, with coefficients
$a_{ij}$ in \eqref{lz} depending continuously on $q,t$.
(They will involve the entries of $\rho_{n}(\sg_{n-1})$ and
the scalar in $A$, which is up to sign a certain power of $q$.)


To show that this linear combination is not a trace multiple
on $\rho_{n-1}(\bt')$ when $\bt'$ ranges over $C$, it suffices,
by continuity, to look at $q=-t=1$. Then the action of
$\sg_i$ is this of permuting the subscripts $i$ and $i+1$ of
the (basis) elements $v_{i,j}$ in the formula \eqref{rf}. Clearly
$\sg_{n-1}$, exchanging subscripts $n-1$ and $n$, does not
fix (or take to multiples of themselves) all such elements
with $1\le i< j\le n-1$. \qed

Theorem \ref{t2} follows by combining the previous three
statements, theorem \ref{t1}, and the result in \cite{reiko}.

\proof[of theorem \reference{t2}]
Let $\bt\in B_{n-1}$ be a braid representation of $L$ as a
non-central braid. Then, by using lemma \ref{l23} (and remark
\ref{z1}), we have that $\rho_{n-1}(\bt)$ is not scalar for proper
$q,t$ of definite form.

If $\psi_{n-1}(\bt)$ is not scalar, then
the claim follows from the work in \cite{reiko}. So assume that
$\psi_{n-1}(\bt)$ is scalar. Then, regarding $B_{n-1}\simeq
B_{1,n-1} \subset B_n$, the map $B_{n-1}\to \bC$ given by
\begin{eqn}\label{ok}
\bt'\mapsto \tr\rho_n(\bt'\sg_{n-1})
\end{eqn}
is for these $q,t$ linear but not a trace multiple by lemma \ref{Q2}.

By theorem \ref{t1}, the closure of the $\rho_{n-1}$-image of the
conjugacy class $C$ of $\bt$ is a $SU(p')$-conjugacy class $D$ with
$p'=(n-1)(n-2)/2$. Since \eqref{ok} is a linear function on $C$,
it can be extended to such a function on $D$, and in a unique way.
This extension is not constant by proposition \ref{prop2}, and $D$
is a connected set. Then this set $D$ cannot contain a dense subset
$C$ on which a continuous map \eqref{ok} takes a finite (or even
discrete) value range.
\qed

Thus we prove in fact a bit more; e.g. for proper $q,t$ the set of
$|\tr\rho_n|$ or $\arg\tr \rho_n$ on $n$-braid representations 
of $L$ has a closure that contains an interval.

\begin{rem}\label{z3}
{}From the perspective of Markov's theorem, it seems more important
to construct \em{irreducible} braids, i.e. such which are \em{not}
conjugate to $\gm\sg_{n-1}^{\pm 1}$ for $\gm\in B_{n-1}$. The examples
in the proof of theorem \reference{t2} can be easily modified by
exchange moves (see \cite{BM}) to ones which at least \em{may}
be potentially irreducible. (Lemma \reference{Q2} needs a
slight adaptation.) But this promises no real advance, as long
as one cannot \em{prove} irreducibility. No decent general
technique exists to establish this property for non-minimal strand
braid representations, except the arguments in \cite{Morton},
which apply in very special cases.
\end{rem}

\noindent{\bf Acknowledgment.}
I would wish to thank to R. Shinjo, S.~Bigelow, R.~Budney, W.~T.~Song,
and I.~Marin 
for some helpful comments, references and discussions. 

{\small

\let\old@bibitem\bibitem
\def\bibitem[#1]{\old@bibitem}


}

\end{document}